\documentclass[10pt]{article}
\usepackage{amsmath}
\usepackage{amssymb}
\usepackage{amsthm}
\usepackage{graphicx}
\usepackage{verbatim}
\usepackage{enumerate}
\font\elevenss=cmss11

\font\eightss=cmss8

\font\sixss=cmss8 at 6pt

\newfam\ssfam
\textfont\ssfam=\elevenss \scriptfont\ssfam=\eightss
 \scriptscriptfont\ssfam=\sixss

\theoremstyle{plain}
\newtheorem{thm}{Theorem}[section]

\newtheorem{lem}[thm]{Lemma}
\newtheorem{pr}[thm]{Proposition}
\newtheorem{cor}[thm]{Corollary}
\newtheorem{defn}[thm]{Definition}

\theoremstyle{remark}

\def\Z{\mathbb{Z}}

\def\R{\mathbb{R}}
\def\C{\mathbb{C}}

\def\S{{\mathcal S}}
\def\Var{{\rm Var} \,}
\def\ee{\varepsilon}
\def\E{{\mathbb E}}
\def\P{{\mathbb P}}
\def\Cox{\hfill \Box}
\def\disp{\displaystyle}
\def\one{{\bf 1}}
\def\|{{\, | \, }}
\def\F{{\mathcal F}}

\def\psit{\tilde{\psi}}
\def\xt{\tilde{x}}
\def\RM{\mathcal{R}_M}
\def\RN{\mathcal{R}_N}
\def\Uni{\mathcal{U}}

\begin{document}
\begin{titlepage}
\begin{center}
{\large \bf Zeros of a random analytic function approach perfect spacing
under repeated differentiation}
\\[5ex]
\end{center}

\begin{flushright}
Robin Pemantle\footnote{Department of Mathematics,
University of Pennsylvania,
209 South 33rd Street,
Philadelphia, PA 19104, USA, {\tt pemantle@math.upenn.edu}}$^,$
\footnote{Research supported in part by NSF grant \# DMS-1209117},
Sneha Subramanian\footnote{Department of Mathematics,
University of California at Irvine,
340 Rowland Hall (Bldg. \# 400),
University of California, 
Irvine, CA 92697-3875,
{\tt ssneha@math.upenn.edu}
}
\end{flushright}

\vfill

\noindent{\bf Abstract:}
We consider an analytic function $f$ whose zero set forms a
unit intensity Poisson process on the real line.  We show that
repeated differentiation causes the zero set to converge in 
distribution to a random translate of the integers. 
\vfill

\\[2ex]


\end{titlepage}

\setcounter{equation}{0}
\section{Introduction}
\label{sec:intro}

Study of the relation of the zero set of a function $f$ to 
the zero set of its derivative has a rich history.  The 
Gauss-Lucas theorem (see, e.g.,~\cite[Theorem~6.1]{marden1949}) 
says that if $f$ is a 
polynomial then the zero set of $f'$ lies in the convex hull 
of the zero set of $f$.  Another property of the differentiation
operator is that it is {\em complex zero decreasing}: the number
of non-real zeros of $f'$ is at most the number of non-real zeros
of $f$.  This property is studied by~\cite{craven-csordas} 
in the more general context of {\em P{\'o}lya-Schur} operators, 
which multiply the coefficients of a power series by a predetermined 
sequence.  Much of the recent interest in such properties of 
the derivative and other operators stem from proposed attacks
on the Riemann Hypothesis involving behavior of zeros under
these operators~\cite{LeMo1974,Conr1983}.  See 
also~\cite[Section~4]{pemantle-hyperbolicity} for a survey 
of combinatorial reasons to study locations of zeros 
such as log-concavity of coefficients~\cite{brenti-LC}
and negative dependence properties~\cite{borcea-branden-liggett}.

The vague statement that differentiation should even out
the spacings of zeros is generally believed, and a number
of proven results bear this out.  For example, a theorem
attributed to Riesz (later rediscovered by others) states
that the minimum distance between zeros of certain entire
functions with only real zeros is increased by differentiation;
see~\cite[Section~2]{farmer-rhoades} for a history of this
result and its appearance in~\cite{Stoy1926} and subsequent
works of J.\ v.\ Sz.--Nagy and of P.\ Walker.

The logical extreme is that repeated differentiation should 
lead to zeros that are as evenly spaced as possible.  If the 
original function $f$ has real zeros, then all derivatives of
$f$ also have all real zeros.  If the zeros of $f$ have some 
long-run density on the real line, then one might expect
the zero set under repeated differentiation to approach a 
lattice with this density.  A sequence of results leading up
to this was proved in~\cite{farmer-rhoades}.  They show that
the gaps between zeros of $f' + af$ are bounded between
the infimum and supremum of gaps between consecutive zeros
of $f$ and generalize this to a local density result that
is applicable to the Riemann zeta function.  They claim a
result~\cite[Theorem~2.4.1]{farmer-rhoades} that implies
the convergence of spacings of zeros to a constant (their
Theorem~2.4.2) but a key piece of their proof, Proposition~5.2.1,
has a hole that seemingly cannot be fixed (D. Farmer, personal 
communication).

The central object of this paper is a random analytic function
$f$ whose zeros form a unit intensity point process.  We 
construct such a function and prove translation invariance 
in Proposition~\ref{th:f and TI}. 
Our main result is that as $k \to \infty$, the zero set of the 
$k^{th}$ derivative of $f$ approaches a random translate of the 
integers.  Thus we provide, for the first time, a proof of
the lattice convergence result in the case of a random zero set.

The remainder of the paper is organized as follows.  
In the next section we give formal constructions and 
statements of the main results.  We also prove preliminary
results concerning the contruction, interpretation and 
properties of the random function $f$.  At the end of
the section we state an estimate on the Taylor coefficients of $f$,
Theorem~\ref{th:a} below, and show that Theorem~\ref{th:main} follows 
from Theorem~\ref{th:a} without too much trouble.  
In Section~\ref{sec:sigma} we begin proving Theorem~\ref{th:a}, that is,
estimating the coefficients of $f$.  It is suggested 
in~\cite{farmer-rhoades} that the Taylor series for $f$ might prove
interesting, and indeed our approach is based on determination of
these coefficients.  We evaluate these via Cauchy's integral formula.
In particular, in Theorem~\ref{th:sigma}, we locate a saddle point 
$\sigma_k$ of $z^{-k} f$.  
In Section~\ref{sec:f} we prove some estimates on $f$, allowing us to 
localize the Cauchy integral to the saddle point and complete the 
proof of Theorem~\ref{th:a}.  
We conclude with a brief discussion.

\setcounter{equation}{0}
\section{Statements and preliminary results}
\label{sec:statements}

We assume there may be readers interested in analytic function
theory but with no background in probability.  We therefore
include a couple of paragraphs of formalism regarding random
functions and Poisson processes, with apologies to those readers 
for whom it is redundant.  

\subsection{Formalities} \label{ss:formal}

A random object $X$ taking values in a $S$ endowed with a 
$\sigma$-field $\S$ is a map $X : (\Omega , \F) \to (S , \S)$ 
where $(\Omega , \F , \P)$ is a probability space.  We will 
never need explicitly to name the $\sigma$-field $\S$ on $\S$, 
nor will we continue to say that maps must be measurable, though 
all maps are assumed to be.  If $S$ is the space of analytic functions, 
the map $X$ may be thought of as a map $f : \Omega \times \C \to \C$.
The statement ``$f$ is a random analytic function'' means that 
for any fixed $\omega \in \Omega$, the function $z \mapsto f(\omega , z)$
is an analytic function.  The argument $\omega$ is always dropped
from the notation, thus, e.g., one may refer to $f'(z)$ or $f(\lambda z)$,
and so forth, which are also random analytic functions.    

A unit intensity Poisson process on the real lines is a random
counting measure $N$ on the measurable subsets of $\R$ such that
for any disjoint collection of sets $\{ A_1 , \ldots , A_n \}$,
of finite measure, the random variables $\{ N(A_1) , \ldots , N(A_n) \}$ 
are a collection of independent Poisson random variables with
respective means $|A_1| , \ldots , |A_n|$ (here $|B|$ denotes the
measure of $B$).  The term ``counting measure'' refers to a measure
taking values in the nonnegative integers; there is a random countable 
set $E$ such that the measure of any set $A$ is the cardinality of 
$A \cap E$.  We informally refer to the set $E := \{ x \in \R : 
N( \{ x \}) = 1 \}$ as the ``points of the Poisson process.''  

Let $\Omega$ henceforth denote the space of counting measures 
on $\R$, equipped with its usual $\sigma$-field $\F$, and let 
$\P$ denote the law of a unit intensity Poisson process.  This
simplifies our notation by allowing us to construct a random 
analytic function $f : \Omega \times \C \to \C$ by a formula
for the value of $f(N , z)$, guaranteeing that the random
function $f$ is determined by the locations of the points of 
the Poisson process $N$.

For $N \in \Omega$ and $\lambda \in \R$, let $\tau_\lambda N$
denote the shift of the measure $N$ that moves points to the right 
by $\lambda$; in other words, $\tau_\lambda N (A) := N (A - \lambda)$ 
where $A - \lambda$ denotes the leftward shift $\{ x - \lambda : x \in A \}$.
A unit intensity point process is translation invariant.  This
means formally that $\P \circ \tau_\lambda = \P$ for any $\lambda$.
If $X$ is a random object in a space $S$ admitting an action 
of the group $(\R , +)$, we say that $X$ is constructed in a
translation invariant manner from $N$ if $X(\tau_\lambda N)
= \lambda X(N)$.  This implies that the law of $X$ is
invariant under the $(\R , +)$-action but not conversely.
In what follows we will construct a random analytic function 
$f$ which is translation invariant up to constant multiple.
Formally, for any function $g$ let $[g]$ denote the set
of functions $\{ \lambda g : g \in \R \}$.  Let $(\R , +)$ act
on the set of analytic functions by translation in the domain:
$\lambda * g(z) := g(z - \lambda)$.  This commutes with the
projection $g \mapsto [g]$.  Our random analytic function $f$ 
will have the property that $[f]$ is constructed in a 
translation invariant manner from $N$.

\subsection{Construction of $f$} 
\label{ss:f}

Various quantities of interest will be defined as sums and
products over the set of points of the Poisson process $N$.
The sum of $g$ evaluated at the points of the counting measure $N$
is more compactly denoted $\int g \, dN$.  If $\int |g| \, dN < \infty$ 
then this is an absolutely convergent sum and its meaning is clear.
Because many of these infinite sums and products are not
absolutely convergent, we introduce notation for some symmetric
sums that are conditionally convergent.

Let $g : \R \to \C$ be any function.  Let $N_M$ denote the restriction 
of $N$ to the interval $[-M,M]$.  Thus, $\int g \, dN_M$ denotes
the sum of $g(x)$ over those points of the process $N$ lying 
in $[-M,M]$.  Define the symmetric integral $\int_* g \, dN$
to be equal to $\lim_{M \to \infty} \int g \, dN_M$ when the
lmit exists.  It is sometimes more intuitive to write such an 
integral as a sum over the points, $x$, of $N$.  Thus we denote
$$\sum_* g(x) := \int_* g(x) \, dN (x) = \lim_{M \to \infty}
   \int g(x) \, dN_M (x)$$
when this limit exists.

Similarly for products, we define the symmetric limit by
$$\prod_* g(s) := 
   \lim_{M \to \infty} \exp \left ( \int \log g \, dN_M \right ) \, .$$
Note that although the logarithm is multi-valued, its integral agains
a counting measure is well definedup to multiples of $2 \pi i$, whence
such an integral has a well defined exponential.

\begin{thm} \label{th:f and TI}
Except for a set of values of $N$ of measure zero, the symmetric product
\begin{equation} \label{eq:f}
f(z) := \prod_* \left ( 1 - \frac{z}{x} \right ) 
\end{equation} 
exists.  The random function $f$ defined by this product is
analytic and translation invariant.  In particular, 
\begin{equation} \label{eq:const}
f(\tau_\lambda N , z) = \frac{f(N , z - \lambda)}{f(N , -\lambda)}
\end{equation}
which implies $[f(\tau_\lambda N , \cdot)] = [f(N , \cdot - \lambda)]$.
\end{thm}

We denote the $k^{th}$ derivative of $f$ by $f^{(k)}$.  The following
is an immediate consequence of Theorem~\ref{th:f and TI}. 
\begin{cor} \label{cor:TI}
For each $k$, the law of the zero set of $f^{(k)} (z)$ is 
translation invariant. 
$\Cox$
\end{cor}

Translation invariance of $f$ is a little awkward because it holds
only up to a constant multiple.  It is more natural to work with
the logarithmic derivative
$$h(z) := \sum_* \frac{1}{z-x} \, .$$

\begin{lem} \label{lem:h TI}
The random function $h$ is meromorphic and its poles are precisely
the points of the process $N$, each being a simple pole.  Also $h$
is translation invariant and is the uniform limit on compact sets 
of the functions 
$$h_M (z) := \int \frac{1}{z-x} \, dN_M (x) \, .$$
\end{lem}

\noindent{\sc Proof:}  
Let $\Delta_M := h_{M+1} (0) - h_M (0)$.  It is easily checked that
\begin{enumerate}[(i)]
\item $\P (\Delta_M > \ee)$ is summable in $M$;
\item $\E \Delta_M = 0$;
\item $\E \Delta_M^2$ is summable.
\end{enumerate}
By Kolmogorov's three series theorem, it follows that $\lim_{M \to \infty}
h_M (0)$ exists almost surely.

To improve this to almost sure uniform convergence on compact sets,
define the $M^{th}$ tail remainder by $T_M (z) := h(z) - h_M (z)$ 
if the symmetric integral $h$ exists.  Equivalently, 
$$T_M (z) := \lim_{R \to \infty} \int \frac{1}{z-x} d(N_R - N_M) (x)$$
if such a  limit exists.  Let $K$ be any compact set of complex numbers.
We claim that the limit exists and that
\begin{equation} \label{eq:claim}
G(M) := \sup_{z \in K} \left | T_M(z) - T_M (0) \right |  
   \to 0 \mbox{ almost surely as } M \to \infty .
\end{equation}
To see this, assume without loss of generality that 
$\disp M \geq 2 \sup \{ |\Re \{ z \}| : z \in K \}$.  Then
\begin{equation} \label{eq:T}
T_M (z) - T_M (0) = \lim_{R \to \infty} \int \left ( 
   \frac{1}{z-x} - \frac{1}{-x} \right ) \, d (N_R - N_M) (x) \, .
\end{equation}
Denote $C_K := \sup_{z \in K} |z|$.  As long as $z \in K$ and 
$|x| \geq M$, the assumption on $M$ gives
\begin{equation} \label{eq:C_K}
\left | \frac{1}{z-x} - \frac{1}{-x} \right | = 
  \left | \frac{z}{x(z-x)} \right | \leq \frac{2C_K}{x^2} \, .
\end{equation}
This implies that the integral in~\eqref{eq:T} is absolutely
integrable with probability~1.  Thus, almost surely,
$T_M (z) - T_M (0)$ is defined by the convergent integral
$$T_M (z) - T_M (0) = \int 
   \left ( \frac{1}{z-x} - \frac{1}{-x} \right ) \, d (N - N_M) (x) \, .$$
Plugging in~\eqref{eq:C_K}, we see that
$\disp G(M) \leq 2 C_K \int x^{-2} \, d(N - N_M) (x) \, ,$ which 
goes to zero (by Lebesgue dominated convergence) except on the
measure zero event that $\int |x|^{-2} \, dN(x) = \infty$.

This proves~\eqref{eq:claim}.  The triangle inequality then yields 
$\sup_{z \in K} |T_M (z)| \leq G(M) + |T_M (0)|$, both summands
going to zero almost surely.  By definition of $T_M$, this means
$h_M \to h$ uniformly on $K$.  The rest is easy.  For fixed $K$
and $M$, $h = h_M + \lim_{R \to \infty} (h_R - h_M)$.  When
$M$ is sufficiently large and $R > M$, the functions $h_R - h_M$
are analytic on $K$.  Thus $h$ is the sum of a meromorphic
function with simple poles at the points of $N$ in $K$
and a uniform limit of analytic functions.  Such a limit
is analytic.  Because $K$ was arbitrary, $h$ is meromorphic 
with simple poles exactly at the points of $N$.

The final conclusion to check is that $h$ is translation invariant.
Unraveling the definitions gives
$$h (\tau_\lambda N , z) = \int_{* \lambda} \frac{1}{(z - \lambda) - x} 
   \, d N (x)$$
where $\int_{* \lambda}$ is the limit as $M \to \infty$ of
the integral over $[-\lambda - M , - \lambda + M]$.  
Translation invariance then follows from checking that
$\int_{M-\lambda}^M \frac{1}{z-x} \, dN(x)$ and
$\int_{-M-\lambda}^{-M} \frac{1}{z-x} \, dN(x)$ both converge
almost surely to zero.  This follows from the large deviation bound 
$$\P \left ( \left | \int_{M-\lambda}^M \frac{1}{z-x} \, dN(x)
   \right | \geq \ee \right ) = O \left ( e^{-c M} \right )$$
and Borel-Cantelli.
$\Cox$

\noindent{\sc Proof of Theorem}~\ref{th:f and TI}:
The antiderivative of the meromorphic function $h$ is an 
equivalence class (under addition of constants) of functions 
taking values in $\C \mod (2 \pi i)$.  Choosing the
antiderivative of $h_m$ to vanish at the origin and exponentiating
gives the functions $f_M$, whose limit as $M \to \infty$ is the 
symmetric product, $f$.  Analyticity follows because $f$ is the
uniform limit of analytic functions.  Translation invariance
up to constant multiple follows from translation invariance of $h$.
The choice of constant~\eqref{eq:const} follows from the 
definition, which forces $f(0) = 1$.
$\Cox$

Before stating our main results, we intoduce a few properties
of the random analytic function $f$.  

\begin{pr} \label{pr:increasing}
$f(\overline{z}) = \overline{f(z)}$ and $|f (a + bi)|$ 
is increasing in $|b|$.
\end{pr}

\noindent{\sc Proof:} Invariance under conjugation is evident
from the construction of $f$.  For $a,b \in \R$,
\begin{eqnarray*}
\log |f (a + bi)| & = & \sum_* \log \left | 1 + \frac{a+bi}{x} \right | \\[2ex]
& = & \frac{1}{2} \sum_* \log \left [ \left ( 1 + \frac{a}{x} \right )^2
   + \left ( \frac{|b|}{x} \right )^2 \right ] \, .
\end{eqnarray*}
Each term of the sum is increasing in $|b|$.
$\Cox$

The random function $f$, being almost surely an entire analytic 
function, almost surely possesses an everywhere convergent 
power series 
$$f(z) = \sum_{n=0}^\infty e_n z^n \, .$$
By construction $f(0) = 1$, hence $e_0 = 1$.  
The function $f$ is the uniform limit on compact sets of
$\disp f_M := \exp \left ( \int \log (1 - z/x) \, dN_M (x) \right )$.
The Taylor coefficients $e_{M,n}$ of $f_M$ are the elementary
symmetric functions of the negative reciprocals of the points
of $N_M$:
$$e_{M,k} = e_k \left ( \{ -1/x : N_M (x) = 1 \} \right ) \, .$$
It follows that $e_{M,k} \to e_k$ as $M \to \infty$ for each fixed $k$.
Thus we may conceive of $e_k$ as the $k^{th}$ elementary symmetric
function of an infinite collection of values, namely the negative 
reciprocals of the points of the Poisson process.  The infinite sum 
defining this symmetric function is not absolutely convergent but
converges conditionally in the manner described above.

We do not know a simple form for the marginal distribution of 
$e_k$ except in the case $k=1$.  To see the distribution of $e_1$,
observe that the negative reciprocals of the points of a unit intensity 
Poisson process are a point process with intensity $dx/x^2$.  
Summing symmetrically in the original points is the same as 
summing the negative reciprocals, excluding those in $[-\ee , \ee]$, 
and letting $\ee \to 0$.
By a well known construction of the stable laws (see, 
e.g.~\cite[Section~3.7]{durrett4}), this immediately implies:

\begin{pr} \label{pr:cauchy}
The law of $e_1$ is a symmetric Cauchy distribution.
$\Cox$
\end{pr}

While we have not before seen a systematic study of 
symmetric functions of points of an infinite Poisson process, 
symmetric functions of IID collections of variables have been 
studied before.  These were first well understood in 
Rademacher variables (plus or minus one with probability $1/2$ each).
It was shown in~\cite[Theorem~1]{mori-szekely} that 
the marginal of $e_k$, suitably normalized, is the value of the 
$k^{th}$ Hermite polynomial on a standard normal random input.
This was extended to other distributions, the most general
result we know of being the one in~\cite{major1999}.  

\subsection{Main result and reduction to coefficient analysis}
\label{ss:reduction}

The random analytic function $f$ is the object of study for the
remainder of the paper.  Our main result is as follows, the proof
of which occupies most of the remainder of the paper.

\begin{thm}[Main result] \label{th:main}
As $k \to \infty$, the zero set of $f^{(k)}$ converges in distribution 
to a uniform random translate of the integers.
\end{thm}

We prove the main result via an analysis of the Taylor coefficients
of $f$, reducing Theorem~\ref{th:main} to the following result.

\begin{thm}[behavior of coefficients of the derivatives] \label{th:a}
Let $a_{k,r} := [z^r] f^{(k)} (z)$.
There are random quantities 
$\{ A_k \}_{k \geq 1}$ and $\{ \theta_k \}_{k \geq 1}$ such that
\begin{eqnarray}
a_{k,r} & = & A_k \left [ \cos \left ( \theta_k - \frac{r \pi}{2} \right )
   + o_k(1) \right ] \, \cdot \frac{\pi^r}{r!} \mbox{ in probability } 
   \label{eq:a1} \\[1ex]
\sum_{r = 1}^\infty M^r \frac{ |a_{k,r}| }{A_k} & < & \infty
   \mbox{ with probability } \, 1 - o(1), \label{eq:a2}
\end{eqnarray}
for any $M > 0$.
The use of the term ``in probability'' in the first statement means 
that for every $\ee > 0$ the quantity
$$\P \left ( \left | \frac{r!}{\pi^r A_k} a_{k,r} 
   - \cos \left( \theta_k - \frac{r \pi}{2} \right)
   \right | > \ee \right )$$
goes to zero for fixed $r$ as $k \to \infty$.
\end{thm}

A surprising consequence of this result is that the signs
of the coefficients $\{ e_k \}$ are periodic with period~4.
In particular, $e_k$ and $e_{k+2}$ have opposite signs
with probability approaching~1 as $k \to \infty$.
It is interesting to compare this with simpler models, 
such as the Rademacher model in~\cite{mori-szekely}
in which a polynomial $g$ has $n$ zeros, each of them 
at $\pm 1$, with signs chosen by independent fair coin flips.  
The number of positive roots will be
some number $b = n/2 + O(\sqrt{n})$.  Once $n$ and $b$
are determined, the polynomial $g$ is equal to $(z-1)^b
(z+1)^{n-b}$.  The coefficients of $g$ are the elementary
symmetric functions of $b$ ones and $n-b$ negative ones.
The signs of these coefficients have 4-periodicity as well
(\cite[Remark~4]{mori-szekely}).  An analogue of Theorem~\ref{th:a} 
in the case of IID variables with a reasonably general common
distribution appears in~\cite{major1999} (see also~\cite{sneha-thesis}
for extensions).  The proofs, in that case as well as in
the present paper, are via analytic combinatorics.  We know
of no elementary argument for the sign reversal between
$e_k$ and $e_{k+2}$.


\noindent{\sc Proof of Theorem~\ref{th:main} from Theorem}~\ref{th:a}:
We assume the conclusion of Theorem~\ref{th:a} holds and establish 
Theorem~\ref{th:main} in the following steps.  Let $\theta_k$ and 
$A_k$ be as in the conclusion of Theorem~\ref{th:a}.

\underline{Step 1: Convergence of the iterated derivatives on compact sets.} 
\footnote{This step is analogous to~\cite[Theorem~2.4.1]{farmer-rhoades},
the correctness of which is unknown to us at this time.}
Let $\psi_k (x) := \cos (\pi x - \theta_k)$.  
Fix any $M > 0$.  Then
\begin{equation} \label{eq:a3}
\sup_{x \in [-M,M]} \left | \frac{f^{(k)} (x)}{A_k} - \psi_k (x) \right |
   \to 0 \mbox{ in probability as } k \to \infty \, .
\end{equation}
To prove this, use the identity $\cos (\theta_k - r \pi / 2)
= (-1)^j \cos (\theta_k)$ when $r = 2j$ and $(-1)^j \sin (\theta_k)$
when $r = 2j+1$ to write
\begin{eqnarray*}
\psi_k (x) & = & \cos(\theta_k) \cos (\pi x) + \sin (\theta_k) \sin (\pi x)
   \\[2ex]
& = & 
   \cos (\theta_k) \left [ 1 - \frac{\pi^2 x^2}{2!} + \cdots \right ] + 
   \sin (\theta_k) \left [ \pi x - \frac{\pi^3 x^3}{3!} + \cdots \right ]  
   \\[2ex]
& = & \sum_{r=0}^\infty \cos \left ( \theta_k - \frac{r\pi}{2} \right )
   \frac{\pi^r}{r!} x^r \, .
\end{eqnarray*}
This last series is uniformly convergent on $[-M,M]$.  Therefore, given
$\ee > 0$ we may choose $L$ large enough so that 
\begin{equation} \label{eq:a1proof}
\sup_{x \in [-M,M]} \left | \psi_k (x) - \sum_{r=0}^L 
   \cos \left ( \theta_k - \frac{r\pi}{2} \right )
   \frac{\pi^r}{r!} x^r \right | < \frac{\ee}{3} \, .
\end{equation}
By ~\eqref{eq:a2}, we may choose $L$ larger if necessary,
in order to ensure that
\begin{equation} \label{eq:a2proof}
\left | \sum_{r=L+1}^\infty \frac{ a_{k,r} }{A_k} x^r \right | < \frac{\ee}{3}
\end{equation}
for all $x \in [-M,M]$.  Fix such an $L$ and use the power series for
$f^{(k)}$ to write
\begin{equation} \label{eq:parts}
\frac{f^{(k)} (x)}{A_r} - \psi_k (x) = 
   \left ( \sum_{r=0}^L \frac{a_{k,r}}{A_k} x^r - \psi_k (x) \right )
   + \sum_{r=L+1}^\infty \frac{a_{k,r}}{A_k} x^r \, .
\end{equation}
Putting~\eqref{eq:a1proof} together with~\eqref{eq:a1} shows that
the first term on the right-hand side of~\eqref{eq:parts} is at
most $\ee / 3 + \sum_{r=0}^L \xi_r$ where $\xi_r$ is the term 
of~\eqref{eq:a1} that is $o_k (1)$ in probability.  By~\eqref{eq:a1}
we may choose $k$ large enough so that 
$\ee / 3 + \sum_{r=0}^L \xi_r < 2 \ee / 3$ with probability at 
least $1 - \ee / 2$.  Thus, we obtain
$$\sup_{x \in [-M,M]} \left | \frac{f^{(k)} (x)}{A_k}
   - \psi_k (x) \right | \leq \ee$$
with probability at least $1 - \ee$, etablishing~\eqref{eq:a3}.

\underline{Step 2: The $k+1^{st}$ derivative as well.} 
Let $\eta_k (x) := - \pi \sin (\pi x - \theta_k)$.  
Fix any $M > 0$.  Then
\begin{equation} \label{eq:a4}
\sup_{x \in [-M,M]} \left | \frac{f^{(k+1)} (x)}{A_k} - \eta_k (x) \right |
   \to 0 \mbox{ in probability as } k \to \infty \, .
\end{equation}
The argument is the same as in Step 1, except that we use the
power series $f^{(k+1)} (x)  =\sum_{r=1}^\infty r a_{k,r} x^{r-1}$
in place of $f^{(k)} (x) = \sum_{r=0}^\infty a_{k,r} x^r$ and
$\disp \eta_k (x) = \sum_{r=1}^\infty \cos (\theta_k - r\pi/2) 
\frac{\pi^r}{(r-1)!} x^{r-1}$.

\underline{Step 3: Convergence of the zero set to some lattice.} 
On any interval $[-M,M]$, the zero set of
$f^{(k)}$ converges to the zero set of $\psi_k$ in probability.
More precisely, for each $\ee > 0$, if $k$ is large enough, then
except on a set of probability at most $\ee$,
for each zero of $f^{(k)}$ in $[-M+2\ee , M-2\ee]$
there is a unique zero of $\psi_k$ within distance $2 \ee$ and 
for each zero of $\psi_k$ in $[-M+2\ee , M-2\ee]$ there is a
unique zero of $f^{(k)}$ within distance $2 \ee$.

This follows from Steps~1 and~2 along with the following fact
applied to $\psi = \psi_k$, $\psit = f^{(k)}$, $I = [-M,M]$
and $c = 1/2$.

\begin{lem} \label{lem:unif zeros}
Let $\psi$ be any function of class $C^1$ on an interval $I := [a,b]$.
Suppose that $\min \{ \|psi| , |\psi'| \} \geq c$ on $I$.  For any
$\ee > 0$, let $I^\ee$ denote $[a+\ee , b - \ee]$.  Let $\ee < c^2$
be positive, and suppose that a $C^1$ function $\psit$ satisfies
$|\psit - \psi| \leq \ee$ and $|\psit' - \psi'| \leq c/2$ on $I$.
Then the zeros of $\psi$ and $\psit$on $I$ are in correspondence
as follows.
\begin{enumerate}[(i)]
\item For every $x \in I^{\ee/c}$ with $\psi (x) = 0$ there is an
$\xt \in I$ such that $\psit (\xt) = 0$ and $|\xt - x| \leq \ee / c$.
This $\xt$ is the unique zero of $\psit$ in the connected component 
of $\{ |\psi| < c \}$ containing $x$.
\item For every $\xt \in I^{\ee/c}$ with $\psit (\xt) = 0$ there is a 
$x \in I$ with $\psi (x) = 0$.  This $x$ is the unique zero of $\psi$
in the connected component of $\{ |\psi| < c \}$ containing $x$.
\end{enumerate}
\end{lem}

\noindent{\sc Proof:} For~$(i)$, pick any $x \in I^{\ee / c}$ with 
$\psi (x) = 0$.  Assume without loss of generality that $\psi' (x) > 0$
(the argument when $\psi' (x) < 0$ is completely analogous).
On the connected component of $|\psi| \leq c$ one has $\psi' > c$.
Consequently, moving to the right from $x$ by at most $\ee / c$
finds a value $x_2$ such that $\psi (x_2) \geq \ee$, 
moving to the left from $x$ by at most $\ee / c$ finds a value 
$x_2$ such that $\psi (x_2) \leq -\ee$, and $\psi'$ will be
at least $c$ on $[x_1 , x_2]$.
We have $|\psit - \psi| \leq \ee$, whence $\psit (x_1) \leq 0 
\leq \psit (x_2)$, and by the Intermediate Value Theorem $\psit$ 
has a zero $\xt$ on $[x_1,x_2]$.  To see uniqueness, note that
if there were two such zeros, then there would be a zero of $\psit'$,
contradicting $|\psit' - \psit| < c/2$ and $|\psi'| > c$.

To prove~$(ii)$. pick $\xt \in I^{\ee / c}$ with $\psit (\xt) = 0$.  
Then $|\psi (\xt)| \leq \ee \leq c$ whence $|\psi' (\xt)| > c$.  
Moving in the direction of decrease of $|\psi (\xt)|$, 
$|\psi'|$ remains at least $c$, so we must hit zero within a
distance of $\ee / c$.  Uniqueness follows again because another
such zero would imply a critical point of $\psi$ in a region 
where $|\psi| < c$.
$\Cox$

\underline{Step 4: Uniformity of the random translation.}  
Because convergence
in distribution is a weak convergence notion, it is equivalent
to convergence on every $[-M,M]$.  We have therefore proved that 
the zero set of $f^{(k)}$ converges in distribution to a random 
translate of the integers.  On the other hand, Corollary~\ref{cor:TI}
showed that the zero set of $f^{(k)}$ is translation invariant for all $k$.
This implies convergence of the zero set of $f^{(k)}$ to a 
uniform random translation of $\Z$, and completes the proof of 
Theorem~\ref{th:main} from Theorem~\ref{th:a}.
$\Cox$

\setcounter{equation}{0}
\section{Estimating coefficients}
\label{sec:sigma}

\subsection{Overview}

The coefficients $e_k := [z^k] f(z)$ will be estimated via the 
Cauchy integral formula
\begin{equation} \label{eq:cauchy}
e_k = \frac{1}{2 \pi i} \int z^{-k} f(z) \frac{dz}{z} \, .
\end{equation}
Denote the logarithm of the integrand by $\phi_k (z) := \log f (z) - 
k \log z$.  Saddle point integration theory requires the 
identification of a saddle point $\sigma_k$ and a contour
of integration $\Gamma$, in this case the circle through $\sigma_k$
centered at the origin, with the following properties.
\begin{enumerate}[$(i)$]
\item $\sigma_k$ is a critical point of $\phi$, that is, 
$\phi' (\sigma_k) = 0$.
\item The contribution to the integral from a arc of $\Gamma$ of
length of order $\phi'' (\sigma_k)^{-1/2}$ centered at $\sigma_k$
is asymptotically equal to 
$e^{\phi (\sigma_k)} \sqrt{2 \pi / \phi''(\sigma_k)}$.
\item The contribution to the integral from the complement of
this arc is negligible.  
\end{enumerate}
In this case we have a real function $f$ with two complex conjugate 
saddle points $\sigma_k$ and $\overline{\sigma_k}$.  Accordingly, there
will be two conjugate arcs contributing two conjugate values to the 
integral while the complement of these two arcs contributes negligibly.
One therefore modifies $(i)$--$(iii)$ to:
\begin{enumerate}[$(i')$]
\item $\sigma_k$ and $\overline{\sigma_k}$ are critical points of $\phi$, 
and there are no others on the circle $\Gamma$, centered at the origin,
of radius $|\sigma_k|$.
\item The contribution to the integral from a arc of $\Gamma$ of
length of order $\phi'' (\sigma_k)^{-1/2}$ centered at $\sigma_k$
is asymptotically equal to 
$e^{\phi (\sigma_k)} \sqrt{2 \pi / \phi''(\sigma_k)}$.  
\item The contribution to the integral from the complement of
the two conjugagte arcs is negligible compared to the contribution 
from either arc.  
\end{enumerate}
Note that $(iii')$ leaves open the possibility that the two contributions
approximately cancel, leaving the supposedly negligible term dominant.

\subsection{Locating the dominant saddle point}
\label{ss:sigma}

The logarithm of the integrand in~\eqref{eq:cauchy}, also known
as the phase function, is well defined up to multiples of $2 \pi i$.  
We denote it by
$$\phi_k (z) := - k \log z + \sum_* \log \left ( 1 - \frac{z}{x} 
   \right ) \, .$$
When $k=0$ we denote $\sum_* \log (1 - z/x)$ simply by $\phi (z)$.

\begin{pr} \label{pr:derivs}
For each $k,r$, the $r^{th}$ derivative $\phi_k^{(r)}$ of the 
phase function $\phi_k$ is the meromorphic function
defined by the almost surely convergent sum
\begin{equation} \label{eq:phi r}
\phi_k^{(r)} (z) = (-1)^{r-1} (r-1)! \left [ - \frac{k}{z^r} 
   + \sum_* \frac{1}{(z-x)^r} \right ] \, .
\end{equation}
Thus in particular,
$$\phi_k' (z) = - \frac{k}{z} + \sum_* \frac{1}{z-x} \, .$$
\end{pr}

\noindent{\sc Proof:} When $r=1$, convergence of~\eqref{eq:phi r}
and the fact that this is the derivative of $\phi$ is just 
Lemma~\ref{lem:h TI} and the subsequent proof of 
Theorem~\ref{th:f and TI} in which $f$ is constructed from $h$.
For $r \geq 2$, with probability~1 the sum is absolutely convergent.
$\Cox$

The main work of this subsection is to prove the following result,
locating the dominant saddle point for the Cauchy integral.

\begin{thm}[location of saddle] \label{th:sigma}
Let $E_{M,k}$ be the event that $\phi_k$ has a unique zero,
call it $\sigma_k$, in the ball of radius $M k^{1/2}$ about 
$i k / \pi$.  Then $\P (E_{M,k}) \to 1$ as $M,k \to \infty$ with $k \geq 
4\pi^{2}M^{2}$.  
\end{thm}

This is proved in several steps.  We first show that 
$\phi_k' (i k / \pi)$ is roughly zero, then use estimates
on the derivatives of $\phi$ and Rouch{\'e}'s Theorem to
bound how far the zero of $\phi_k'$ can be from $i k / \pi$.

The function $\phi_k'$ may be better understood if one applies the
natural scale change $z = ky$.  Under this change of variables,
$$\phi_k' (z) = - \frac{1}{y} + \sum_* \frac{1/k}{y - x/k} \, .$$ 
Denote the second of the two terms by
$$h_k (y) := \sum_* \frac{1/k}{y - x/k} \, .$$ 
We may rewrite this as $\disp h_k (Y) = \int \frac{1}{y - x} dN^{(k)} (x)$
when $N^{(k)}$ denotes the rescaled measure defined by 
$N^{(k)} (A) = k^{-1} N(kA)$.  The points of the process $N^{(k)}$
are $k$ times as dense and $1/k$ times the mass of the points of $N$.
Almost surely as $k \to \infty$ the measure $N^{(k)}$ converges to
Lebesgue measure.  In light of this it is not surprising that 
$h_k (y)$ is found near $\disp \int \frac{1}{z-y} \, dy$. 
We begin by rigorously confirming this, the integral being
equal to $- \pi \, {\rm sgn} \, \Im \{ z \}$ away from the real axis.

\begin{lem} \label{lem:interchange}
    If $z$ is not real then
    \begin{equation*}
	\E \int_{*} \frac{1}{|z - x|^{m}} dN(x) = \lim_{M\to\infty} \E \int 
	\frac{1}{|z - x|^{m}} dN_{M}(x),
    \end{equation*}
    for $m \geq 2$, and
    \begin{equation*}
	\E \int_{*} \frac{1}{(z - x)^{m}} dN(x) = \lim_{M\to\infty} \E \int 
	\frac{1}{(z - x)^{m}} dN_{M}(x),
    \end{equation*}
    for $m \geq 1$.
\end{lem}

\noindent{\sc Proof:}
The first equality holds trivially by Monotone Convergence Theorem. Next, write 
$\RM$ as the number of points of the process $N$ withing $[-M,M]$, and $L = 
2\Im(z)$. Then, for $m = 2$,
\begin{eqnarray*}
    \E \int \frac{1}{|z - x|^{2}} dN_{M}(x) & = & \E \sum_{j: |X_{j}| \leq M} 
    \frac{1}{\Im(z)^{2} + (\Re(z) - X_{j})^{2}} \\[1ex]
    & \leq & \E \left( \frac{\mathcal{R}_{L}}{\Im(z)^{2}} \right) + \, \E 
    \sum_{j: L \leq |X_{j}| \leq M} \frac{4}{|X_{j}|^{2}} \\[1ex]
    & \leq & \frac{2L}{\Im(z)^{2}} + \frac{4}{L}.
\end{eqnarray*}
Therefore, as $\Im(z) \neq 0$, $\E \int_{*} \frac{1}{|z - x|^{2}} dN(x) < 
\infty$, and moreover, $\E \int_{*} \frac{1}{|z - x|^{m}} dN(x) < \infty, 
\forall m \geq 2$. Thus, by Dominated Convergence Theorem,
 \begin{equation*}
     \E \int_{*} \frac{1}{(z - x)^{m}} dN(x) = \lim_{M\to\infty} \E \int 
	\frac{1}{(z - x)^{m}} dN_{M}(x)
\end{equation*}
holds for $m \geq 2$. We shall now show the above to hold true for $m = 1$.

Note that
\begin{equation*}
    \E \left[ \left| \int \frac{1}{z - x} dN_{M}(x) \right|^{2} \right] =  
    \E \int \frac{1}{|z - x|^{2}} dN_{M}(x) \, + \, \E \sum_{j \neq k: |X_{j}|, 
    |X_{k}| \leq M} \frac{1}{(z - X_{j})(\overline{z} - X_{k})}.
\end{equation*}
The first term in the above equation converges to $\E \int_{*} \frac{1}{|z - 
x|^{2}} dN(x)$ as $M \to \infty$. As for the second part,
\begin{equation*}
    \E \sum_{j \neq k: |X_{j}|, |X_{k}| \leq M} \frac{1}{(z - 
    X_{j})(\overline{z} - X_{k})} = \E \left[ \RM (\RM - 1) \cdot \E \left( 
    \frac{1}{(z - \Uni_{2})(\overline{z} - \Uni_{2})} \right) 
    \right],
\end{equation*}
where $\Uni_{1}$ and $\Uni_{2}$ are i.i.d. Uniform$(-M,M)$ random variables. So,
\begin{eqnarray*}
    \E \sum_{j \neq k: |X_{j}|, |X_{k}| \leq M} \frac{1}{(z - 
    X_{j})(\overline{z} - X_{k})} & = & \left( \int_{-M}^{M} \frac{1}{z - u} du 
    \right)^{2} \\[1ex]
    & = & \left[ - \log \left| \frac{M - z}{M + z} \right| -
    i \arctan \left ( \frac{M - \Re(z)}{\Im(z)} \right) \right. \\
    & & \left. + \, i \arctan \left( \frac{- M - \Re(z)}{\Im(z)} \right) 
    \right]^{2} \\[1ex]
    & \longrightarrow & - \pi^{2}, \text{ as } M \to \infty.
\end{eqnarray*}
Thus the quantities $\left \{ \E \left[ \left| \int \frac{1}{z - x} dN_{M}(x) 
\right|^{2} \right], M > 0 \right \}$ have a uniform upper bound - let us call 
it $B(z)$. Then, given $\ee > 0$,
\begin{eqnarray*}
    \E \left[ \left| \int \frac{1}{z - x} dN_{M}(x) \right| \cdot \one_{\left| 
    \int \frac{1}{z - x} dN_{M}(x) \right| \geq K}\right] & \leq & \frac{1}{K} 
    \cdot \E \E \left[ \left| \int \frac{1}{z - x} dN_{M}(x) \right|^{2} \right]
     \\[1ex]
    & \leq & \frac{B(z)}{K} < \ee,
\end{eqnarray*}
for $K > \frac{B(z)}{\ee}$. Thus, if $z$ is not real, $\left \{ \E \left[ 
 \int \frac{1}{z - x} dN_{M}(x) \right], M > 0 \right \}$ is a uniformly 
 integrable collection, and hence, converges in $L_{1}$.
$\Cox$

\begin{pr} \label{pr:1}
If $z$ is not real then
\begin{equation} \label{eq:1}
\E \left [ \int_* \frac{1}{z - x} \, dN(x) \right ] = \mp i \pi
\end{equation}
with the negative sign if $z$ is in the upper half plane and the
positive sign if $z$ is in the lower half plane.
If $z$ is not real and $m \geq 2$ then
\begin{equation} \label{eq:2}
\E \left [ \int_* \frac{1}{(z - x)^m} \, dN(x) \right ] = 0 \, .
\end{equation}
\end{pr}

\noindent{\sc Proof:} If $\RM$ denotes the number of Poisson points in [-M,M], 
then conditioning on $\RM$, the poisson points $X_j$ that are contained in 
$[-M,M]$ are identically and independently distributed as Uniform$[-M,M]$. Then,
\begin{equation*}
\E \left[ \left. \int \frac{1}{z - x} \, dN_M(x) \right| \RM 
\right] = \RM \cdot \E \left ( \frac{1}{z-\Uni} \right ),
\end{equation*}
where $\Uni \sim$ Uniform$[-M,M]$. Writing $z = r e^{i\theta}$, we get,
\begin{eqnarray*}
\E \left[ \left. \int \frac{1}{z - x} \, dN_M(x) \right| \RM 
\right] & = & \frac{\RM}{2M} \int_{x \in [-M,M]} \frac{1}{r\cos\theta + i 
r\sin\theta - x}dx\\[1ex]
& = & \RM\left[ \frac{-1}{2M} \log \left| \frac{M - z}{M + z} \right| -
\frac{i}{2M} \arctan \left ( \frac{M - r\cos\theta}{r\sin\theta} \right ) 
\right.\\
& &  + \, \left. \frac{i}{2M} \arctan \left( \frac{- M - 
r\cos\theta}{r\sin\theta} \right) \right]\\[1ex]
\implies \E \left[ \int \frac{1}{z - x} \, dN_M(x) \right] & = &  - 
\log \left| \frac{M - z}{M + z} \right| - i \arctan \left ( \frac{M -
r\cos\theta}{r\sin\theta} \right )\\
& &  + \, i \arctan \left ( \frac{- M - r\cos\theta}{r\sin\theta} \right )
\end{eqnarray*}
since, $\RM \sim $ Poisson$(2M)$. Taking $M \to \infty$, by Lemma 
\ref{lem:interchange} we get,
\begin{equation*}
\E \left [ \int_* \frac{1}{z - x} \, dN(x) \right ] = - \pi i,
\end{equation*}
for $z$ in the upper half plane, and,
\begin{equation*}
\E \left [ \int_* \frac{1}{z - x} \, dN(x) \right ] = \pi i,
\end{equation*}
for $z$ in the lower half plane, where the interchange of limits and 
expectation is by Lemma \ref{lem:interchange}.

Now fix $m \geq 2$ and $z \notin \R$.
\begin{eqnarray*}
\E \left[ \left. \int \frac{1}{(z - x)^m} \, dN_M(x) \right| \RM \right]
& = & \RM \cdot \E \left[ \frac{1}{(z - \Uni)^m} \right]\\[1ex]
& = & \frac{\RM}{2M} \cdot \frac{1}{m - 1} \left\{ \frac{1}{(z-M)^{m-1}} -
\frac{1}{(M+z)^{m-1}} \right\}.\\[1ex]
\implies \E \left[ \int \frac{1}{(z - x)^m} \, dN_M(x) \right] & = &
\frac{1}{m-1}\left\{ \frac{1}{(z-M)^{m-1}} - \frac{1}{(M+z)^{m-1}} \right\}.
\end{eqnarray*}
Thus, using Lemma \ref{lem:interchange},
\begin{equation*}
\E \left[ \int_* \frac{1}{(z - x)^m} \, dN(x) \right] = \lim_{M \to \infty}
\frac{1}{m-1}\left\{ \frac{1}{(z-M)^{m-1}} - \frac{1}{(M+z)^{m-1}} \right\} = 0.
\end{equation*}
$\Cox$

The next proposition and its corollaries help us to control how 
much the functions $\phi_k$ and $h_k$ can vary.  These will be used 
first in Lemma~\ref{lem:4.3.6}, bounding $h_k$ over a ball, then
in Section~\ref{ss:saddle est} to estimate Taylor series
involving $\Phi_k$.  We begin with a general result on the variance
of a Poisson integral.

\begin{pr} \label{pr:3.5a}
Let $\psi : \R \to \C$ be any bounded function with $\int |\psi (x)|^2
\, dx < \infty$.  Let $Z$ denote the compensated Poisson integral
of $\psi$, namely
$$Z := \lim_{M \to \infty} \left [ 
   \int \psi (x) \, dN_M (x) - \int_{-M}^M \psi (x) \, dx \right ] \, .$$
Then $Z$ is well defined and has finite variance given by
$$\E |Z|^2 = \int |\psi (x)|^2 \, dx \, .$$
\end{pr}

\noindent{\sc Proof:} This is a standeard result but the proof
is short so we supply it.  Let
$$Z_M := \int \psi (x) \, dN_M (x) - \int_{-M}^M \psi (x) \, dx $$
and let $\Delta_M := Z_M - Z_{M-1}$ denote the increments.  We apply
Kolmogorov's Three Series Theorem to the independent sum 
$\sum_{M=1}^\infty \Delta_M$, just as in the proof of Lemma~\ref{lem:h TI}.  
Hypothesis (i) is satisfied because $\int_M^{M+1} |\psi|$ goes to zero.  
Hypothesis (ii) is satisfied because $\E \Delta_M = 0$ for all $M$.
To see that hypothesis (iii) is satisfied, observe that $\E |\Delta_M|^2
= \int |\psi (x)|^2 \one_{M-1 \leq |x| \leq M} \, dx$, the summability of
which is equivalent to our assumption that $\psi \in L^2$.  
We conclude that the limit exists almost surely.  By monotone 
convergence as $M \to \infty$, $\Var (Z) = \int |\psi|^2$.
$\Cox$

Define 
$$W_r (z) := \int_* (z-x)^{-r} \, dN(x) \, .$$
If $\alpha > 1$ and $\lambda$ is real, the intergal 
$\int |z-x|^{-\alpha} \, dx$ is invariant under $z \mapsto z + \lambda$ 
and scales by $\lambda^{1-\alpha}$ under $z \mapsto \lambda z$.  
Plugging in $\psi (x) = (z - x)^{-r}$ therefore yields the following 
immediate corollary.  
\begin{cor} \label{cor:3.5b}
Let $z$ have nonzero imaginary part and let $r \geq 2$ be an integer.
Then $W_r (z)$ is well defined and there is a positive constant $\gamma_r$
such that
$$\E |W_r (z)|^2 = \frac{\gamma_r}{|\Im \{ z \}|^{2r-1}} .$$
$\Cox$
\end{cor}

In the case of $r=1$ we obtain the explicit constant $\gamma_1 = 1$:
\begin{equation*}
\E |W_1(z) \mp \pi i |^2 = \frac{\pi}{|\Im(z)|}.
\end{equation*}
To see this, compute
\begin{eqnarray*}
\E \left[ \left. \int \frac{1}{| z - x |^2} \, dN_M(x) \right| \RN \right] & = & \RN
\cdot \frac{1}{2N} \int_{x \in [- N, N]} \frac{1}{(z - x) \cdot (\bar{z} - x)}
dx\\[1ex]
& = & \frac{\RN}{2N (\bar{z} - z)} \left [ \int_{x \in [- N, N]} \frac{1}{z -
x}dx - \int_{x \in [- N, N]} \frac{1}{\bar{z} - x}dx \right ]\\[1ex]
& = & \frac{1}{\bar{z} - z} \left \{ \E \left [ \left. \int \frac{1}{z - x} \, dN_M(x)
\right| \RN \right ] - \E \left [ \left. \int \frac{1}{\bar{z} - x} \, dN_M(x) \right|
\RN \right ] \right \}.
\end{eqnarray*}
Thus, taking expectations and by Lemma \ref{lem:interchange}
\begin{equation*}
\E \left [ \int_* \frac{1}{|z - x|^2} \, dN(x) \right ] = \frac{1}{\bar{z} - z} \left
\{ \E \left [ \int_* \frac{1}{z - x} \, dN(x) \right ] - \E \left [ \int_*
\frac{1}{\bar{z} - x} dN(x) \right ]  \right \}.
\end{equation*}
Proposition \ref{pr:1} shows the difference of expectations on the
right-hand side to be $- 2 i \pi$, yielding $\gamma_1 = \pi$.

\begin{cor} \label{cor:3.5c}
For $y$ with nonzero imaginary part and $r \geq 1$, $W_r (ky)$
has variance $\E [ \Re \{ W - \overline{W} \}^2 + 
\Im \{ W - \overline{W} \}^2] = k^{-1/2} \gamma_r (y)$.  
It follows (with $\delta_{1,r}$ denoting the Kronecker delta), that
\begin{equation} \label{eq:W_r}
\phi_k^{(r)} (ky) = - i \pi \delta_{1,r} + (r-1)! k^{1-r} 
   \left ( \frac{-1}{y} \right )^r + O \left ( k^{1/2 - r} \right )
\end{equation}
in probability as $k \to \infty$. 
\end{cor}

\noindent{\sc Proof:}  Let $N^{(k)}$ denote a Poisson law of intensity
$k$, rescaled by $k^{-1}$.  In other words, $N^{(k)}$ is the average
of $k$ independent Poisson laws of unit intensity.  Under the change 
of variables $u = x/k$, the Poisson law $dN (x)$ becomes $k dN^{(k)} (u)$.
Therefore,
\begin{eqnarray*}
W_r (ky) & = & \int_* \frac{1}{(ky-x)^r} \, dN(x) \\[2ex]
& = & k^{1-r} \int_* \frac{1}{(y-u)^r} \, dN^{(k)} (u) \\[2ex]
& = & k^{1-r} \left ( \frac{1}{k} \sum_{j=1}^k W_r^{[j]} \right )
\end{eqnarray*}
where $\{ W_r^[1] , \ldots , W_r^{[k]} \}$ are $k$ independent
copies of $W_r (y)$.  Because $W_r (y)$ has mean $- i \pi \delta_{1,r}$
and variance $\gamma_r (y)$, the variance of the
average is $k^{-1/2} \gamma_r (y)$.  The remaining conclusion
follows from the expression~\eqref{eq:phi r} for $\phi_k^{(r)}$ and
the fact that a random variable with mean zero and variance
$V$ is $O(V^{1/2})$ in probability. 
$\Cox$

\subsection{Uniformizing the estimates}

At some juncture, our pointwise estimates need to be strengthened 
to uniform estimates.  The following result is a foundation for 
this part of the program.

\begin{lem} \label{lem:lipschitz}
Fix a compact set $K$ in the upper half plane and an integer $r \geq 1$.
There is a constant $C$ such that for all integers $k \geq 1$,
$$\E \sup_{z \in K} \left | h_k^{(r)} (z) \right | \leq C k^{-1/2} \, .$$
\end{lem}

\noindent{\sc Proof:} Let $F^{(k)}$ denote the CDF for the 
random compensated measure $N^{(k)} - dx$ on $\R^+$, thus 
$F(x) = N^{(k)}[0,x] - x$ when $x > 0$ and $F(x) = x - N^{(k)}[x,0]$ 
when $x < 0$.  We have 
$$h_k^{(r)} (z) = \int_* C (z-x)^{-r-1} \, dN(x) = \int_* C (z-x)^{-r-1}
   \, dF^{(k)} (x)$$
because $\int_* (z-x)^{-r-1} \, dx = 0$.  This leads to 
\begin{eqnarray*}
&& \E \sup_{z \in K} \left | h_k^{(r)} (z) \right |  \\[1ex]
   & \leq & \lim_{M \to \infty}
   \E \sup_{z \in K} \left | \int_0^M \frac{1}{(z-x)^r} 
      \, dF^{(k)}(x) \right | + 
   \E \sup_{z \in K} \left | \int_{-M}^0 \frac{1}{(z-x)^r} \, dF^{(k)}(x) 
      \right | \, .
\end{eqnarray*}
The two terms are handled the same way.  Integrating by parts, 
$$\int_0^M (z-x)^{-r} \, dF^{(k)}(x) = (z-x)^{-r} N[0,M] - 
   \int - r (z-x)^{-r-1} \, F^{(k)} (x) \, dx .$$
This implies that 
\begin{eqnarray*}
&& \E \sup_{z \in K} \left | h_k^{(r)} (z) \right | \\[2ex]  
& \leq & \lim_{M \to \infty} \left [ \E |F^{(k)}(M)| \sup_{z \in K} |z-x|^{-r} 
   + \int_0^M  \sup_{z \in K} r |z-x|^{-r-1} |F^{(k)} (x)| 
   \right ] \, dx \\[2ex]
& \leq & C_K \lim_{M \to \infty}
   \left ( M^{-r + 1/2} + k^{-1/2} \right ) \, .
\end{eqnarray*}
Sending $M$ to infinity gives the conclusion of the lemma.
$\Cox$

\begin{cor} \label{cor:lipschitz} ~~\\[-5ex]
\begin{enumerate}[(i)]
\item $\disp \sup_{z \in K}|h_k^{(r)} (z)| = O(k^{-1/2})$ in probability.
\item $h_k$ and its derivatives are Lipschitz on $K$ with Lipschitz
constant $O(k^{-1/2})$ in probability.
\item For $r \geq 2$, the $O(k^{-1/2})$ term in the 
expression~\eqref{eq:W_r} for $\phi_k^{(r)} (ky)$ is uniform 
as $y$ varies over compact sets of the upper half plane.
\end{enumerate}
\end{cor}

\noindent{\sc Proof:} Conclusion~$(i)$ is Markov's inequality. 
Conclusion~$(ii)$ follows because any upper bound on a function
$|g'|$ is is a Lipschitz constant for $g$.  Conclusion~$(iii)$ 
follows from the relation between $h_k$ and $\phi_k$.  
$\Cox$

\begin{lem} \label{lem:4.3.6}
For any $c > 0$,
$$\P \left [ \sup \left \{ \left |h_k (y) + i \pi \right| \, : \, |y - 
\frac{i}{\pi}| \leq M k^{-1/2} \right \} \geq c M k^{-1/2} \right ] \to 0$$
as $M \to \infty$ with uniformly in $k \geq 4 \pi^2 M^2$.
\end{lem}

\noindent{\sc Proof:} Fix $c , \ee > 0$.  Choose $L > 0$ such that
the probability of the event $G$ is at most $\ee/2$, where $G$ is 
the event that the Lipschitz constant for some $h_k$ on the ball 
$B(i\pi , 1/(2\pi))$ is greater than $L$.  Let $B$ be the ball
of radius $M k^{-1/2}$ about $i/\pi$.  The assumption $k \geq
4 \pi^2 M^2$ guarantees that $B$ is a subset of the ball 
$B(i\pi , 1 / (2\pi))$ over which the Lipschitz constant 
was computed.  Let $y$ be any point in $B$. 
The ball of radius $\rho := c M k^{-1/2} \ee / (2L)$ about $y$ 
intersects $B$ in a set whose area is at least $\rho^2 \sqrt{3}/2$, 
the latter being the area of two equilateral triangles of side $\rho$.  
If $|h_k (y) + i / \pi| \geq c M k^{-1/2}$ and $G$ goes not occur, 
then $|h_k (u) + i / \pi| \geq (1/2) c M k^{-1/2}$ on the ball of 
radius $\rho$ centered at $y$.  

Now we compute in two ways the expected measure $\E |S|$ of the set 
$S$ of points $u \in B$ such that $|h_k (u) + i \pi| \geq \rho$.
First, by what we have just argued, 
\begin{equation} \label{eq:Q1}
\E |S| \geq \frac{\sqrt{3}}{2} \rho^2 \left ( Q - \frac{\ee}{2} \right )
   = \left ( Q - \frac{\ee}{2} \right ) \sqrt{\frac{3 c^2 \ee^2}{16 L^2}}
   \frac{M^2}{k}
\end{equation}
where $Q$ is the probability that there exists a $y \in B$ such that
$|h_k (y) + i / \pi| \geq c M k^{-1/2}$.  Secondly, 
by Proposition~\ref{pr:1} and the computation of $\gamma_1$, 
for each $u \in B$, $\E h_k (u) + i/\pi = 0$ and 
$\E |h_k (u)|^2 = \pi / k$, leading to 
$\E |h_k (u) + i/\pi| \leq \sqrt{2 \pi / k}$ and hence
\begin{eqnarray*}
\P \left ( \left | h_k (u) + \frac{i}{\pi} \right | \geq \rho \right )
   & \leq & \frac{\sqrt{2 \pi/k}}{\rho} \\
   & = & \frac{\sqrt{2 \pi / k}}{c M k^{-1/2} \ee / (2L)} \\[1ex]
   & = & \sqrt{\frac{32 \pi L^2}{c^2}} M^{-1/2} \, .
\end{eqnarray*} 
By Fubini's theorem, 
\begin{equation} \label{eq:Q2}
\E |S| \leq |B| \sup_{u \in B} \P \left ( 
   \left | h_k (u) + \frac{i}{\pi} \right | \geq c M k^{-1/2} \right ) 
   \leq \pi \frac{M^2}{k} \sqrt{\frac{32 \pi L^2}{c^2}} M^{-1/2} \, .
\end{equation}
Putting together the inequalities~\eqref{eq:Q1} and~\eqref{eq:Q2} gives
$$Q - \frac{\ee}{2} \leq \sqrt{\frac{512 \pi^3 L^4}{3 c^4 \ee^2}} 
   M^{-1/2} \, .$$
Once $M$ is sufficiently larger that the radical is at most $\ee / 2$,
this implies that $Q \leq \ee$.  Because $\ee > 0$ was arbitrary,
we have shown that $Q \to 0$ as $M \to \infty$ uniformly in $k$,
as desired.
$\Cox$

\noindent{\sc Proof of Theorem}~\ref{th:sigma}:
Using Lemma \ref{lem:4.3.6} for $c < 1$, we know that
\begin{equation*}
\P \left [ \sup \left \{ |h_k (y) + i \pi| \, : \, |y - \frac{i}{\pi}| \leq M 
  k^{-1/2} \right \} \leq c M k^{- 1/2} \right ] \longrightarrow 1, \text{ as, }
  k\to\infty.
\end{equation*}
Writing
$$ A_{M,k} = \left \{ \omega : \sup \left \{ |h_k (y) + i \pi| \, : \, |y - 
  \frac{i}{\pi}| \leq M k^{-1/2} \right \} \leq c M k^{- 1/2} \right \}, $$
$\forall \omega \in A_{M,k}$, and all $y$ such that $|y - \frac{i}{\pi}| =
Mk^{-1/2}$,
\begin{eqnarray*}
\left| \phi_k(y) (\omega) - \left( - i\pi - \frac{1}{y} \right) \right| & = &
|h_k(y)(\omega) + i\pi| \\[1ex]
& \leq & cMk^{-1/2} \\[1ex]
& = & c \left | y - \frac{i}{\pi} \right | \\[1ex]
& < & \left | y - \frac{i}{\pi} \right | 
\end{eqnarray*}
for $k$ sufficiently large. Thus, by Rouche's theorem, $ \phi_k(y) (\omega) $ 
and $y - \frac{i}{\pi}$ have the same number of zeros inside 
the disc centered at $i/\pi$ of radius $Mk^{- 1/2}$, i.e. exactly one. This 
implies that, $\P (E_{M,k}) \to 1$ as $M,k \to \infty$ with $k \geq 
4\pi^{2}M^{2}$.
$\Cox$

\setcounter{equation}{0}
\section{The Cauchy integral}

\subsection{Dominant arc: saddle point estimate}
\label{ss:saddle est}

We sum up those facts from the foregoing subsection that 
we will use to estimate the Cauchy integral in the dominant 
arc near $\sigma_k$.

\begin{lem} \label{lem:summarizing phi}
~~\\[-4ex]
\begin{enumerate}[(i)]
\item $\phi' (\sigma_k) = 0$.
\item $\sigma_k^2 \phi'' (\sigma_k) = k + O(k^{1/2})$ in probability
as $k \to \infty$.
\item If $K$ is the set $\{ z : |z - \sigma_k| \leq k/2$ then
$\disp \sup_{z \in K} k^3 \phi^{(3)} (z) = O(k)$ in probability.
\end{enumerate}
\end{lem}

\noindent{\sc Proof:} 
The first is just the definition of $\sigma_k$.  For the second, 
using Corollary~\ref{cor:lipschitz} for $r = 2$ and 
$y = \frac{i}{\pi}$, the estimate~\eqref{eq:W_r} is uniform, hence
$$\phi'' (\sigma_k) = \phi'' \left ( \frac{ik}{\pi} \right )
   + O \left ( k^{-3/2} \right ) = \frac{-\pi^2}{k} 
   + O \left ( k^{-3/2} \right )$$
in probability.  Multiplying by $\sigma_k^2 \sim -k^2 / \pi^2$ 
gives~$(ii)$.  The argument for part~$(iii)$ is analogous to the
argument for part~$(ii)$.
$\Cox$

\begin{defn}[Arcs, fixed value of $\delta$]
For the remainder of the paper, fix a number $\delta \in (1/3 , 1/2)$.
Parametrize the circle $\Gamma$ through $\sigma_k$ in several pieces,
all oriented counterclockwise, as follows (see Figure~\ref{fig:Gamma}).  
Define $\Gamma_1$ to be 
the arc $\{z: z = \sigma_{k}e^{it}, - k^{-\delta} \leq t \leq k^{-\delta} \}$.
Define $\Gamma_1'$ to be the arc $\{z: z = \overline{\sigma}_{k}e^{it}, 
- k^{-\delta} \leq t \leq k^{-\delta} \}$, so that the arc is conjugate 
to $\Gamma_1$ but the orientation remains  counterclockwise.  Define 
$\Gamma_2$ to be the part of $\Gamma$ in the second quadrant that
is not part of $\Gamma_1$, define $\Gamma_3$ to be the part of
$\Gamma$ in the first quadrant not in $\Gamma_1$, and define 
$\Gamma_2'$ and $\Gamma_3'$ to be the respective conjugates.
Define the phase function along $\Gamma$ by 
$$g_k (t) := \phi_{k} (\sigma_k e^{it}) \, .$$
\begin{figure}[!ht] \label{fig:Gamma}
\centering
\includegraphics[width=2.8in]{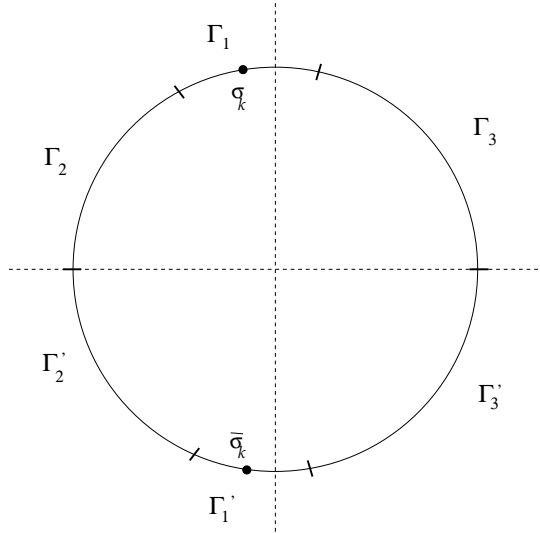}
\caption{Parametrization of the circular contour $\Gamma$}
\end{figure}
\end{defn}

\begin{thm}[Contribution from $\Gamma_1$] 
\label{th:good arcs}
For any integer $r \geq 0$,
$$\frac{\disp{\int_{\Gamma_1} \frac{f(z)}{z^{k + r + 1}} dz}}
   {\disp{k^{-1/2} f(\sigma_k) \sigma_k^{-k-r}}} 
   \longrightarrow i \sqrt{2 \pi}$$
in probability as $k \to \infty$.
\end{thm}

\noindent{\sc Proof:} 
For fixed $k$, the Taylor's expansion of $g_{k}(t)$ gives us,
\begin{equation*}
g_k(t) = g_{k}(0) + t g_k'(0) + \frac{t^2}{2} g_k^{(2)}(0) + \frac{t^3}{6}
\left( \Re g_{k}^{(3)}(t_1) + i \Im g_{k}^{(3)}(t_2)\right),
\end{equation*}
where $t_1$ and $t_2$ are points that lie between $0$ and $t$.\\
By Lemma \ref{lem:summarizing phi}, $g_{k}'(0) = 0$ and
$$g_k^{2}(0) = k + O\left ( k^{1/2} \right )$$ 
in probability.  Thus,
\begin{equation*}
\sup_{|t| \leq k^{- \delta}} \sqrt{k} \left[ \exp \left( \frac{t^2}{2} 
g_k^{(2)}(0) \right) - \exp \left( - \frac{kt^2}{2} \right) \right] 
\longrightarrow 0.
\end{equation*}
In addition, 
\begin{equation*}
\sup_{z \in \Gamma_1} \left| \frac{\sigma_k^r}{z^r} - 1 \right| \longrightarrow
0,
\end{equation*}
while, Lemma \ref{lem:summarizing phi} also gives us
\begin{equation*}
\sup_{|t| \leq k^{- \delta}} \frac{t^3}{6} g_k^{(3)}(t) \longrightarrow 0.
\end{equation*}
Thus,
$$\int_{\Gamma_1} \frac{f(z)}{z^{k + r + 1}} dz \; = \; 
   i \int_{- k^{- \delta}}^{k^{- \delta}} \sigma_k^{- r}
   \exp \left[ g_{k}(0) + \frac{t^2}{2} g_k^{(2)}(0) + \frac{t^3}{6} 
   \left( \Re g_{k}^{(3)}(t_1) + i \Im g_{k,N}^{(3)}(t_2)\right) - i \, r t \right] dt,
   $$
whence, as $k \to \infty$,
$$\sqrt{k} \frac{\int_{\Gamma_1} 
   \frac{f(z)}{z^{k + r + 1}} dz}{\sigma_k^{- r} \exp(g_{k}(0))} 
   - i \sqrt{k} \int_{- k^{- \delta}}^{k^{- \delta}} 
   \exp \left(- \frac{kt^2}{2} \right) dt  \longrightarrow 0.  $$
Changing variables to $t = u/\sqrt{k}$ shows that when $\delta < 1/2$,
the integral is asymptotic to $\sqrt{2 \pi/k}$.  Plugging in 
$g_k (0) = f(\sigma_k) \sigma_k^{- k}$ completes the proof.
$\Cox$

\subsection{Negligble arcs and remainder of proof of 
\protect{Theorem~\ref{th:a}}}
\label{sec:f}

We now show that the Cauchy integral receives negligible contributions
from $\Gamma_2, \Gamma_2', \Gamma_3$ and $\Gamma_3'$.  By conjugate
symmetry we need only check $\Gamma_2$ and $\Gamma_3$; the arguments
are identical so we present only the one for $\Gamma_2$.  

Let $R := |\sigma_k|$ and let $\beta$ denote the polar argument
of $\sigma_k$, that is, $\beta := \arg (\sigma_k) - \pi / 2$,
so that $\sigma_k = i R e^{i \beta}$.  By Theorem~\ref{th:sigma}, 
$\beta = O(k^{-1/2})$ in probability.  We define an exceptional
event $G_k$ of probability going to zero as follows.
\begin{quote}
Let $G_k$ be the event that either $R \notin [k/(2\pi) , 2 k/\pi]$
or $\beta > k^{-\delta} / 2$. 
\end{quote}
If $z = i R e^{i \theta}$ is a point of $\Gamma_2$ with polar argument 
$\theta$, then $\theta$ is at least $k^{-\delta} - |\beta|$, hence 
is at least $(1/2) k^{-\delta}$ on $G_k^c$.
Note that the notation suppresses the dependence of $R$ and $\beta$ 
on $k$, which does not affect the proof of the in-probability result 
in Lemma~\ref{lem:bad arcs}.

\begin{lem} \label{lem:bad arcs}
\begin{equation} \label{eq:big}
\frac{\disp{\int_{\Gamma_2} \frac{f(z)}{z^{k + r + 1}} dz}}
   {\disp{k^{-1/2} f(\sigma_k) \sigma_k^{-k-r}}} 
   \longrightarrow 0
\end{equation}
in probability as $k \to \infty$.
\end{lem}

\noindent{\sc Proof:} Let $z = i R e^{i \theta} \in \Gamma_2$.
Our purpose is to show that $|f(z) z^{-k}|$ is much smaller
that $|f (\sigma_k) \sigma_k^{-k}|$.  On $\Gamma_2$ we are
worried only about the magnitude, not the argument, so we
may ignore the $z^{-k}$ and $\sigma^{-k}$ terms, working 
with $\phi$ rather than with $\phi_k$.  
This simplifies~\eqref{eq:W_r} to
\begin{equation} \label{eq:Re}
\phi' (z) = - i \pi + O \left ( k^{-1/2} \right ) \,
\end{equation}
the estimate being uniform on the part of $\Gamma_2$ with
polar argument less than $\pi/2 - \ee$ by part~$(iii)$ of
Corollary~\ref{cor:lipschitz}.  Let $H_k$ be the exceptional
event where the constant in the uniform $O(k^{-1/2})$ term 
is greater than $k^{1/2-\delta} / 100$, the probability of 
$H_k$ going to zero according to the corollary.  

Integrating the derivative of 
$\Re \{ \phi (z) \}$ along $\Gamma$ then gives
\begin{equation} \label{eq:logs}
\log \left | \frac{f(z)}{f(\sigma_k)} \right | = 
   \pi \left ( \Im \{ z \} - \Im \{ \sigma_k \} \right )
   + O \left ( k^{-1/2} |z - \sigma_k| \right ) \, .
\end{equation}
The first of the two terms is $\pi R (\cos (\theta) - \cos (\beta))$
which is bounded from above by $- (R/2) (\theta^2 - \beta^2)$
which is at most $-(R/4) \theta^2$ on $G_k^c$.  The second term 
is at most 
$$\frac{k^{1/2 - \delta}}{100} k^{-1/2} (2 R \theta)$$ 
on $G_k^c \cap H_k^c$, provided that $\theta \leq \pi/2 - \ee$.
Combining yields 
$$\log \left | \frac{f(z)}{f(\sigma_k)} \right | \leq 
   - \frac{R}{4} \theta^2 + \frac{k^{-\delta}}{100} (2 R \theta) \, \leq \, 
   - R \theta (\left ( \frac{\theta}{4} - \frac{k^{-\delta}}{50} \right ) 
   \, \leq \, - \frac{R \theta^2}{8}$$
on $\Gamma_2$ as long as the polar argument of $z$ is at most 
$\pi/2 - \ee$.  Decompose $\Gamma_2 = \Gamma_{2,1} + \Gamma_{2,2}$ 
according to whether $\theta$ is less than or greater than $\pi/2 - \ee$.
On $G_k^c$ we know that $\theta \geq (1/2) k^{-\delta}$ and 
$R \geq k/(2\pi)$, hence on $\Gamma_{2,1}$,
$$\log \left | \frac{f(z)}{f(\sigma_k)} \right | \leq 
   - \frac{k^{1 - 2\delta}}{64 \pi} \, .$$

Using $d\theta = dz / z$ we bound the desired integral from above by 
$$\left | \frac{\disp{\int_{\Gamma_2} \frac{f(z)}{z^{k + r + 1}} dz}}
   {\disp{k^{-1/2} f(\sigma_k) \sigma_k^{-k-r}}} \right |
   \leq
   \sqrt{k} \int_{\Gamma_2} \left | \frac{f(z)}{f(\sigma_k)} \right |
   \, d\theta \, .$$
On $G_k^c \cap H_k^c$, the contribution from $\Gamma_{2,1}$ is at most
\begin{equation} \label{eq:Gamma_1}
\sqrt{k} \, |\Gamma_2|  \, \exp \left [ - \frac{k^{1 - 2\delta}}{64 \pi} 
   \right ] \, .
\end{equation}

Finally, to bound the contribution from $\Gamma_{2,2}$, use 
Proposition~\ref{pr:increasing} to deduce $|f(z)| \leq |f(z')|$
where $\Re \{ z' \} = \Re \{ z \}$ and $\Im \{ z' \} = k/(4 \pi)$.
Integrating~\eqref{eq:Re} on the line segment between $\sigma_k$ 
and $z'$ now gives~\eqref{eq:logs} again, and on $G_k^c \cap H_k^c$
the right-hand side is at most $- (k/4) + k^{-\delta} k  <  -k/8$ 
once $k \geq 8$.  This shows the contribution from $\Gamma_{2,2}$
to be at most $\ee R e^{-k/8}$.  Adding this to~\eqref{eq:Gamma_1}
and noting that $\P (G_c \cup H_k) \to 0$ proves the lemma.
$\Cox$

\begin{thm} \label{th:e_k}
For fixed $r$ as $k \to \infty$,
$$e_{k+r} = 2 (- 1)^{k + r} \, \Re \left \{ (1 + o(1)) \sigma_k^{-k-r} \, 
f(\sigma_k) \, \sqrt{\frac{1}{2 \pi k}} \right \}$$
in probability as $k \to \infty$.
\end{thm}

\noindent{\sc Proof:} By Cauchy's integral theorem,
$$e_{k+r} = \frac{(- 1)^{k + r}}{2 \pi i} \int_{\Gamma}
   f(z) z^{-k-r-1} \, dz \, .$$
By Theorem~\ref{th:good arcs} and the fact that the contributions 
from $\Gamma_1$ and $\Gamma_1'$ are conjugate, their sum is
twice the real part of a quantity asymptotic to
\begin{equation} \label{eq:leading}
\frac{1}{\sqrt{2 \pi k}} f(\sigma_k) \sigma_k^{-k-r} \, .
\end{equation}
By Lemma~\ref{lem:bad arcs}, the contributions from the remaining
four arcs are negligible compared to~\eqref{eq:leading}.  
The theorem follows.
$\Cox$

\noindent{\sc Proof of Theorem}~\ref{th:a}: 
By the definition of $a_{k,r}$, using Theorem~\ref{th:e_k} to
evaluate $e_k$,
\begin{eqnarray*}
a_{k, r} & = & (- 1)^{k + r} \, e_{k + r} \, \frac{(k+r)!}{r!} \\[1ex]
& = & 2 \, k! \, \frac{(k+r)!}{k!} \, \frac{1}{r!} \, \Re \left \{ (1 + o(1))  
  \sigma_k^{-k-r} \, f(\sigma_k) \, \sqrt{\frac{1}{2 \pi k}} \right \} \, .
\end{eqnarray*}
For fixed $r$ as $k \to \infty$ asymptotically $(k+r)! / k! \sim k^r$.
Setting $\disp A_k = k! \sqrt{\frac{2}{\pi k}} \left| \sigma_k^{-k}
f(\sigma_k) \right|$ and $\disp \theta_k = \arg \{ \sigma_k^{- k} 
f(\sigma_k) \}$ simplifies this to
$$ A_{k} \frac{k^r}{|\sigma_k|^r} 
   \left [ \cos \left ( \theta_k - r \arg (\sigma_k) \right ) \right ] \, .$$
Because in probabiltiy $\arg (\sigma_k) = \pi / 2 + o(1)$ while
$|\sigma_k| \sim k / \pi$, this simplifies finally to
$$a_{k,r} = A_k \left [ \cos \left ( \theta_k - \frac{r \pi}{2} \right )
   + o(1) \right ] \cdot \frac{\pi^r}{r!} \mbox{ in probability},$$
proving the first part of the theorem.

Next, from the proof of Theorem \ref{th:good arcs} it is clear that
\begin{equation*}
\left| \frac{ \int_{\Gamma_1} \frac{f(z)}{z^{k + r + 1}} dz}
{\frac{f(\sigma_k)}{\sigma_k^{k + r}}} \right| \leq \int_{\Gamma_1} 
| \exp(g_k(t) - g_k(0)) | dt
\end{equation*}
is bounded above in probability, and this bound is independent of $r$.
Also the convergence in the proof of Lemma \ref{lem:bad arcs} is 
independent of $r$. Therefore,
\begin{equation*}
\left| \frac{a_{k,r}}{A_k} \right| = O \left( \frac{(k+r)!}{k!}  
\frac{1}{r! |\sigma_k|^r} \right).
\end{equation*}
Since, for any $M > 0$,
\begin{equation*}
\sum_{r=1}^\infty \frac{(k+r)!}{k!} \frac{\pi^r}{r!} \frac{M^r}{k^r} < \infty, \,
\forall \, k > M\pi,
\end{equation*}
with the convergence being uniform over $k \in [T, \infty)$, with $T > M\pi$,
we have our result.

$\Cox$

\clearpage
\bibliographystyle{alpha}
\bibliography{RP}

\end{document}